\documentclass[a4paper,fleqn]{cas-sc}

\usepackage{setspace}
\usepackage{bm}
\usepackage{multirow}

\usepackage[numbers,sort&compress]{natbib}
\usepackage{graphicx}
\usepackage{geometry}
\usepackage{enumitem}
\usepackage{amsmath}
\usepackage{times}
\usepackage{caption}
\usepackage[ruled]{algorithm2e}
\usepackage[title]{appendix}%

\def\tsc#1{\csdef{#1}{\textsc{\lowercase{#1}}\xspace}}
\tsc{WGM}
\tsc{QE}
\tsc{EP}
\tsc{PMS}
\tsc{BEC}
\tsc{DE}
\begin{document}

\let\WriteBookmarks\relax
\def\floatpagepagefraction{1}
\def\textpagefraction{.001}
\shorttitle{Design of Teardrop Hovering Formation along the NRHO}
\shortauthors{Shuyue Fu et~al.}

\title [mode = title]{Design and Continuation of Nonlinear Teardrop Hovering Formation along the Near Rectilinear Halo Orbit}                      

\author[1,2]{Shuyue Fu}[orcid=0009-0001-5111-9779,style=chinese]
\ead{fushuyue@buaa.edu.cn}
\credit{Data curation, Formal analysis, Methodology, Software, Writing - Original draft preparation, Writing - review $\&$ editing}

\author[1,2]{Yihan Peng}[style=chinese]
\ead{pengyihan@buaa.edu.cn}
\credit{Data curation, Formal analysis, Writing - review $\&$ editing}

\affiliation[1]{organization={Shen Yuan Honors College, Beihang University},
                addressline={Xueyuan Road No.37}, 
                city={Beijing},
                postcode={100191}, 
                country={People's Republic of China}}

\affiliation[2]{organization={School of Astronautics, Beihang University},
                addressline={Xueyuan Road No.37}, 
                postcode={100191}, 
                city={Beijing},
                country={People's Republic of China}}

\author[2]{Shengping Gong}[style=chinese]
\cormark[1]
\ead{gongsp@buaa.edu.cn}
\credit{Conceptualization, Funding acquisition, Writing - review $\&$ editing}

\author[2]{Peng Shi}[style=chinese]
\ead{shipeng@buaa.edu.cn}

\credit{Methodology, Formal analysis, Writing - review $\&$ editing}

\cortext[cor1]{Corresponding author}

\begin{abstract}
This short communication is devoted to the design and continuation of a teardrop hovering formation along the Near Rectilinear Halo orbit and provides further insights into future on-orbit services in the cislunar space. First, we extend the concept of the teardrop hovering formation to scenarios along the Near Rectilinear Halo orbit in the Earth-Moon circular restricted three-body problem. Then, we develop two methods for designing these formations based on the nonlinear model for relative motion. The first method addresses the design of the teardrop hovering formations with relatively short revisit distances, while the second method continues hovering trajectories from short to longer revisit distances. In particular, new continuation method is developed to meet the design requirements of this new scenario. Simulation results verify the effectiveness of the proposed methods, and a near-natural teardrop hovering formation is achieved by considering the dynamical properties near the NRHO. Comparisons between design results obtained using linear and nonlinear models further strengthen the necessity of using the nonlinear model.
\end{abstract}

\begin{highlights}
\item The teardrop hovering formation along the NRHO is investigated.
\item New continuation method is developed to meet the design requirements of this scenario.
\item Dynamical properties near the NRHO is used to achieve the near-natural teardrop hovering formation.
\end{highlights}


\maketitle

\section{Introduction} \label{sec1}
As several recent lunar missions (e.g., \textit{Artemis} (USA) \cite{smith2020artemis} and \textit{Chang’e} (China) \cite{Zheng2023}) have been proposed and executed, there has been a renewed interest in lunar exploration that has further promoted the resource utilization of the Moon and cislunar space. The \textit{Gateway}, an international lunar space station and an important component of the \textit{Artemis} program, offers opportunities for scientific research and serves as a crucial staging point for future Mars exploration \cite{smith2020artemis}. As a subset of the Halo orbit families in the Earth-Moon circular restricted three-body problem (CR3BP), Near Rectilinear Halo orbits (NRHOs) are characterized by stability or near-stability, making them suitable for long-term lunar missions performed in the real multi-body environment \cite{zimovan2019dynamical,sanchez2020chance}. For this consideration, the NRHO has been selected as the baseline orbit of the \textit{Gateway} \cite{lee2019white}. To support the practical mission, research on the NRHO has attracted significant attention. Recent studies on the NRHO have mainly focused on its dynamical theory \cite{zimovan2019dynamical,villegas2023resonant}, the design of transfer trajectories associated with the NRHO \cite{trofimov2020transfers,ueda2021multi}, stationkeeping in the multi-body environment \cite{guzzetti2017stationkeeping,davis2017orbit}, and relative motion along the NRHO \cite{sanchez2020chance,sandel2024natural}. Among these topics, the studies on the relative motion along the NRHO provide the theoretical foundation for on-orbit service of the spacecraft moving in the NRHO \cite{sandel2024natural}. On-orbit service, including rendezvous, proximity operation, docking, and hovering, has been thoroughly studied for low Earth orbits in the context of the two-body dynamics. However, as human exploration of the Moon continues, an increasing number of missions are considering the deployment of spacecraft in periodic orbits in the cislunar space (e.g., the aforementioned NRHO), which makes on-orbit service technologies in the cislunar space increasingly important. Despite this, on-orbit service along the NRHO in the CR3BP, particularly hovering formations, remains less explored. Therefore, this short communication focuses on the design of the hovering formation along the NRHO. 

The hovering formation refers to a scenario in which a spacecraft performs control inputs to maintain a specific distance relative to a target \cite{bai2020teardrop,bai2022practical}, such as asteroids \cite{lee2014almost,wang2020asteroid} and spacecraft moving in the reference orbits (i.e., the chief spacecraft) \cite{zhang2013characteristic,dang2014modeling}. Among these configurations, the teardrop hovering formation, named for its teardrop-like trajectory, features hovering trajectories that revisit a fixed position relative to the chief spacecraft after a specific time interval (i.e., the revisit period). To achieve this revisit, only an impulse should be performed during a single revisit period \cite{bai2022practical}. This type of formation not only ensures a stable flyaround of the chief spacecraft but also facilitates rendezvous operations \cite{bai2020teardrop,sandel2024natural}. For the teardrop hovering formations along two-body orbits, Zhang et al. \cite{zhang2016flyaround} designed the impulsive teardrop hovering trajectories based on the analytical solution of the Clohessy-Wiltshire (CW) equations, while Prince and Cobb \cite{prince2018optimal} employed the teardrop relative trajectories to achieve quasi-hover around resident space objects moving in the geostationary orbit. The aforementioned works focused on teardrop formations along circular orbits; moreover, Bai et al. \cite{bai2020teardrop,bai2022practical} developed methods to design teardrop hovering formations along elliptic orbits under J2-perturbed dynamics. To support future on-orbit service missions along the NRHO, we extend the concept of the teardrop hovering formation to the NRHO scenarios. To our best knowledge, this scenario has not been thoroughly investigated. Recent studies on relative motion and on-orbit service along the NRHO have mainly focused on rendezvous \cite{sanchez2020chance,colagrossi2021guidance} and loitering \cite{sandel2024natural,scheuerle2025relative}, whose relative trajectory configurations differ from those of the teardrop hovering formation. Moreover, due to the eigenvalue distribution of the monodromy matrix (i.e., state transition matrix after one orbital period) of the NRHO \cite{fossa2022two}, there is a potential to achieve a natural or near-natural 1:1 teardrop hovering formation (i.e., formation whose revisit period is set equal to the period of the NRHO) with very low fuel consumption by considering the near-periodical properties of the orbit \cite{franz2022database}. Therefore, we focus on the 1:1 teardrop hovering formation in this short communication.

Designing hovering trajectories along the NRHO necessitates consideration of the models for relative motion in the CR3BP. Differing from the relative motion in the two-body dynamics, the relative motion in the CR3BP lacks an analytical solution and should be studied numerically. Furthermore, the models for relative motion in the CR3BP can be categorized into linear model \cite{sanchez2020chance} and nonlinear model \cite{khoury2020orbital}. In regions near the perilune of the NRHO, the linear model may yield significant errors due to noticeable nonlinear effects, high velocities, and close proximity to the Moon \cite{jenson2024bounding,cunningham2023interpolated}. Consequently, we adopt the nonlinear model for relative motion to design the teardrop hovering formation along the NRHO. Two numerical methods to design the teardrop hovering formation are proposed: one for designing teardrop hovering formations with relatively short revisit distances, and another for those with relatively long revisit distances. In the first method, the linear model is used to generate initial guesses for the teardrop hovering formation, followed by differential correction to satisfy the revisit constraints in the nonlinear model. In the second method, based on teardrop hovering formation obtained from the first method, a continuation method is developed to obtain the hovering trajectories from short to long revisit distances. Typically, the numerical continuation method is a predictor-corrector method \cite{capdevila2016a}. In the first step of the continuation method, the feasible direction is predicted based on the existing feasible solution. Then, the initial guess of the next feasible solution is based on prediction and corrected to satisfy the constraints in the second step. Previous works on the continuation method applied to the multi-body problems mainly focused on the continuation of periodic orbits \cite{oshima2022continuation,fu2025analytical}, constrained transfer trajectories \cite{topputo2013optimal,oshima2019low} in initial states and parameters, and continuation of the constrained trajectories from the Earth-Moon planar circular restricted three-body problem to other high-fidelity models \cite{singh2021low}. However, differing from the continuation scenarios reported in the aforementioned works, the continuation scenario requiring to address in this short communication is the continuation in a parameter that is part of the constraints, i.e., the revisit distance. For a similar scenario (continuing transfer trajectories in altitude and inclination of the target orbits), Singh et al. \cite{singh2019mission} employed the natural parameter continuation method. This method makes minor adjustments to the initial states or parameters based on existing feasible solutions and obtains other feasible solutions through differential correction. However, it does not provide an analytical prediction of the feasible direction and may suffer from relatively poor convergence performance and only a partial mapping of the solution space \cite{fu2025analytical}. In contrast, the linear predictor method can theoretically provide the continuation direction and yield a more accurate prediction of the feasible direction \cite{van2016tadpole}. Therefore, we derive a new form of the linear predictor using the least square method to address this new continuation scenario. Simulation results verify the proposed methods, and the comparison between design results obtained from the linear and nonlinear models confirms the necessity of using the nonlinear model to design teardrop hovering formations along the NRHO. The main contributions of this short communication can be summarized as follows:

\begin{enumerate}[label=(\arabic*)]
\item We extend the concept of the teardrop hovering formation to NRHO scenarios and develop two methods for designing teardrop hovering formations along the NRHO. In particular, we propose a new continuation method to meet the new design requirements of this scenario.
\item The impulse distribution in our design results reflects the dynamical properties near the NRHO, thereby providing valuable insights into parameter selection for the design of future on-orbit service missions along the NRHO. By considering these dynamical properties, a near-natural teardrop hovering formation is achieved.
\item  We adopt the nonlinear model for relative motion in the CR3BP to design the teardrop hovering formation along the NRHO, and simulations confirm the necessity of using the nonlinear model.
\end{enumerate}

The rest of this short communication is organized as follows. Section \ref{sec2} presents the the mathematical background, including the CR3BP, NRHO, and the models for relative motion in the CR3BP. Section \ref{sec3} introduces the concept of the teardrop hovering formation and develops two methods for designing hovering trajectories. The proposed methods are verified in Section \ref{sec4}. Finally, conclusions are drawn in Section \ref{sec5}.

\section{Mathematical Background}\label{sec2}
In this section, we present the mathematical background of this short communication, including the CR3BP model, NRHO, and linear and nonlinear models for relative motion in the CR3BP.
\subsection{CR3BP and NRHO}\label{subsec2.1}
\begin{figure}[h]
\centering
\includegraphics[width=0.27\textwidth]{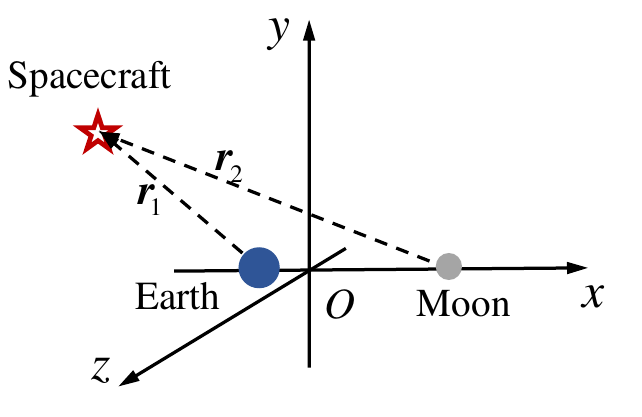}
\caption{The Earth-Moon CR3BP and Earth-Moon rotating frame.}
\label{fig1}
\end{figure}

In this study, we adopt the Earth–Moon CR3BP, where the Earth and Moon are assumed to move in circular orbits around their barycenter. The spacecraft is treated as a massless particle influenced solely by the gravitational forces of the Earth and Moon. The Earth-Moon rotating frame is adopted to describe the dynamical equations, as shown in Fig. \ref{fig1}. 

We consider the case of spatial problem, i.e., the spacecraft can move both within and out of the Earth–Moon orbital plane. The dimensionless units are defined as follows: the length unit (LU) is the Earth-Moon distance, the mass unit (MU) is the combined mass of the Earth and Moon, and the time unit (TU) is ${\text{TU}}={T_\text{EM}}/{2\pi}$, where ${T_{{\text{EM}}}}$ is the period of the motion of the Earth and Moon around their barycenter. Then, the dynamical equations for the spacecraft in the Earth-Moon CR3BP can be written as:
\begin{equation}
\left[ {\begin{array}{*{20}{c}}
  {\begin{array}{*{20}{c}}
  {\dot x} \\ 
  {\dot y} \\ 
  {\dot z} \\
  {\dot u} \\ 
  {\dot v} \\ 
  {\dot w} 
\end{array}} 
\end{array}} \right] = \left[ {\begin{array}{*{20}{c}}
  {\begin{array}{*{20}{c}}
  u \\ 
  v \\ 
  w \\
  {x + 2v - \frac{{\left( {1 - \mu } \right)\left( {x + \mu } \right)}}{{{r_1}^3}} - \frac{{\mu \left( {x + \mu  - 1} \right)}}{{{r_2}^3}}} \\ 
  {y - 2u - \frac{{\left( {1 - \mu } \right)y}}{{{r_1}^3}} - \frac{{\mu y}}{{{r_2}^3}}} \\ 
  { - \frac{{\left( {1 - \mu } \right)z}}{{{r_1}^3}} - \frac{{\mu z}}{{{r_2}^3}}} 
\end{array}} 
\end{array}} \right]\label{eq1}
\end{equation}

\begin{equation}
{r_1} = \sqrt {{{\left( {x + \mu } \right)}^2} + {y^2} + {z^2}} \label{eq2}
\end{equation}

\begin{equation}
{r_2} = \sqrt {{{\left( {x + \mu -1} \right)}^2} + {y^2} + {z^2}} \label{eq3}
\end{equation}
where $\bm{X} = \left[ x,\text{ }y,\text{ }z,\text{ }u,\text{ } v,\text{ }w\right]^{\text{T}}$ is the orbital state, and $\mu $ is the mass parameter expressed as $\mu  = {m_{\text{M}}}/\left( {{m_{\text{E}}} + {m_{\text{M}}}} \right)$. The parameters ${m_{\text{E}}}$ and ${m_{\text{M}}}$ denote the masses of the Earth and Moon, respectively. For numerical trajectory calculations in the Earth–Moon CR3BP, we employ MATLAB®’s ode113 command with the variable step-size, variable order (VSVO) Adams–Bashforth–Moulton algorithm, using absolute and relative tolerances of $1 \times 10^{-13}$. The specific values of the parameters used in the simulations can be found in Table \ref{tab1} \cite{topputo2013optimal}. Additionally, the Earth-Moon CR3BP an invariant known as the Jacobi constant, which is expressed as:
\begin{equation}
C = -\left( {{u^2} + {v^2}}+w^2 \right) +  \left( {{x^2} + {y^2}} \right) + \frac{{2(1 - \mu) }}{{{r_1}}} + \frac{2\mu }{{{r_2}}} + \mu \left( {1 - \mu } \right) \label{eq4}
\end{equation}

\begin{table}[!htb]
\caption{Parameter Setting for the Earth-Moon PCR3BP}\label{tab1}%
\centering
\renewcommand{\arraystretch}{1.5}
\begin{tabular}{@{}llll@{}}
\hline
Symbol & Value  & Units & Meaning\\
\hline
$\mu$    & $1.21506683 \times {10^{ - 2}}$   & --  & Earth-Moon mass parameter  \\
${T_{{\text{EM}}}}$    & $2.24735067 \times {10^6}$   & s  & Earth-Moon period  \\
$R_{\text{E}}$    & $6378.145$   & km  & Mean Earth’s radius  \\
$R_{\text{M}}$    & $1737.100$   & km  & Mean Moon’s radius  \\
LU    & $3.84405000 \times {10^5}$   & km  & Length unit  \\
TU    & $3.75676968 \times {10^5}$   & s  & Time unit  \\
\hline
\end{tabular}
\end{table}

\begin{figure}[h]
\centering
\includegraphics[width=0.6\textwidth]{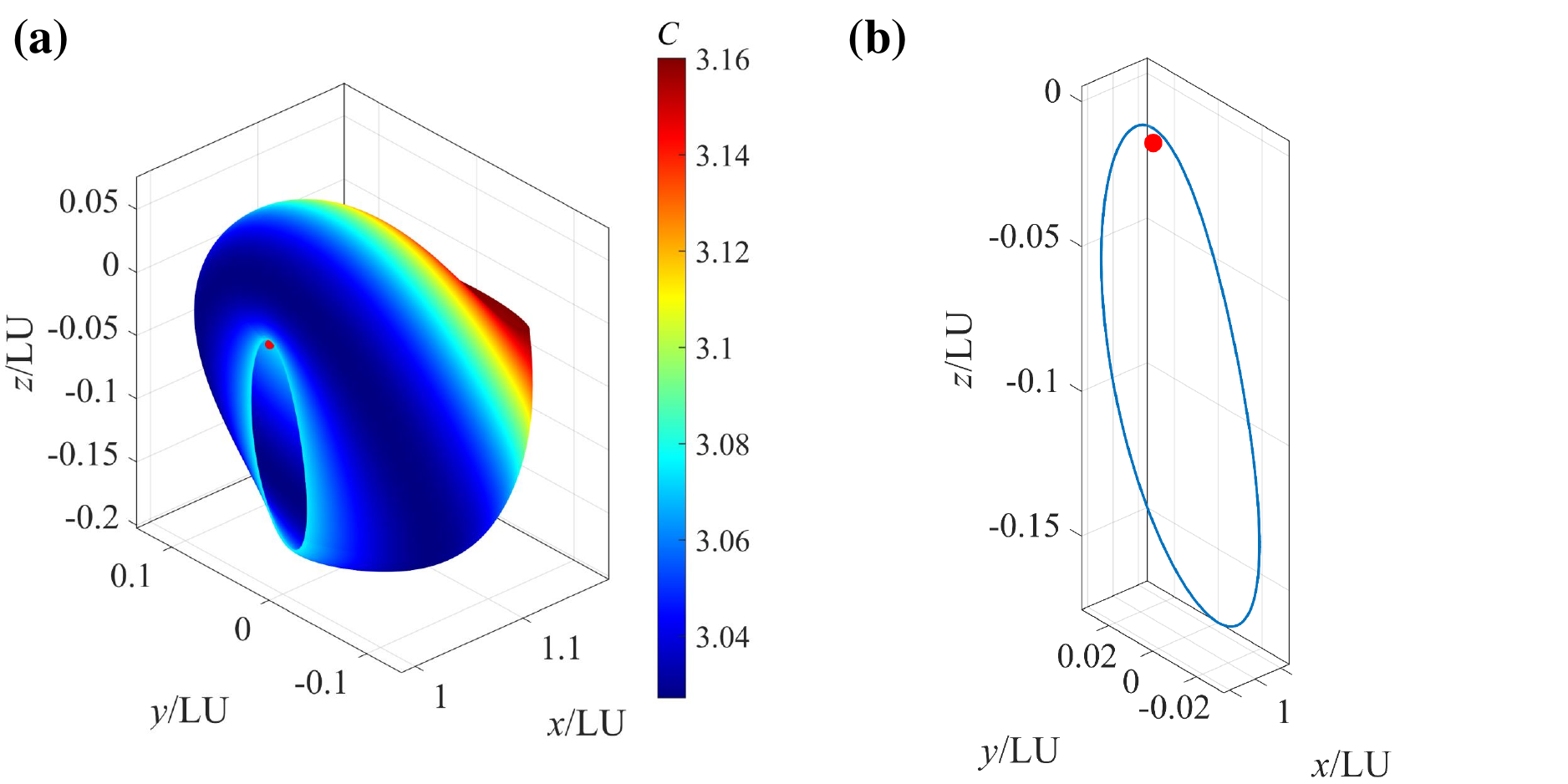}
\caption{The southern L2 family of the Halo orbits and 9:2 NRHO. (a) The southern L2 family of the Halo orbits; (b) 9:2 NRHO.}
\label{fig2}
\end{figure}

In the Earth-Moon CR3BP, several periodic orbits have been discovered, including the libration point periodic orbits, and quasi-satellite orbits \cite{franz2022database}. Among these periodic orbits, the Near Rectilinear Halo Orbit (NRHO), a subset of the Halo orbit families, is selected as the baseline orbit for practical missions like the \textit{Gateway} (USA) \cite{lee2019white}. The NRHO is particularly suitable for long-term lunar exploration due to its stability or near-stability. Moreover, the NRHO associated with the southern L2 family of the Halo orbits offers good coverage for the lunar South pole, supporting scientific observations \cite{sanchez2020chance}. Figure \ref{fig2} (a) shows the southern L2 family of the Halo orbits generated from the linear continuation method mentioned in Ref. \cite{fu2025analytical}. This short communication focuses on the relative motion along the 9:2 NRHO associated with the southern L2 family of the Halo orbits shown in Fig. \ref{fig2} (b), which is the baseline orbit of \textit{Gateway} \cite{lee2019white}. 

The NRHO is an axisymmetric orbit. The initial states and period of the considered NRHO is summarized in Table \ref{tab2} ($T_{\text{NRHO}}$ denotes the period of the NRHO, the subscript “\textit{i}” denotes the initial states, and the subscript “NRHO” denotes the states associated with the NRHO). 

\begin{table}
\caption{Initial States And Period of The Considered NRHO}\label{tab2}%
\centering
\renewcommand{\arraystretch}{1.5}
\begin{tabular}{@{}lll@{}}
\hline
Parameter & Value  & Unit \\
\hline
$x_{i\text{NRHO}}$ & $0.987581435006489$ & LU \\
$z_{i\text{NRHO}}$ & $0.005276210630165$ & LU \\
$v_{i\text{NRHO}}$ & $2.120240531159090$ & LU/TU \\
$T_{\text{NRHO}}$ & $4\pi/9$ & TU \\
\hline
\end{tabular}
\end{table}

The eigenvalues of the monodromy matrix of the considered NRHO exhibit the following properties: $\lambda_1=1/\lambda_2$, ${\lambda _5} = {\bar \lambda _6}$, and $\lambda_3=\lambda_4=1$ \cite{fossa2022two}. Among these eigenvalues, $\lambda_3=\lambda_4=1$ is associated with the periodical properties of the NRHO, implying the existence of periodic or near-periodic orbits near the NRHO \cite{franz2022database}. This dynamical property facilitates the design of a natural or near-natural 1:1 teardrop hovering formation along the NRHO. To design teardrop hovering formations along the considered NRHO, the models for relative motion in the CR3BP are first introduced, including the linear model and nonlinear model, which are described in the variational framework.

\subsection{Linear Model for Relative Motion in the CR3BP}\label{subsec2.2}
Several models have been proposed to describe the relative motion along a reference orbit in the CR3BP, including the linear model \cite{chihang2023close}, the second-order nonlinear model \cite{cunningham2023interpolated}, and the full nonlinear model \cite{sandel2024natural}. This short communication focuses on teardrop hovering formations constructed in the full nonlinear model. However, since the full nonlinear model lacks an analytical solution, numerical methods such as differential correction and numerical continuation are used to design the hovering trajectories. The linear model, in contrast, can provide an initial guess to be corrected to the full nonlinear model. Therefore, the linear model is first introduced. 

For the relative motion analysis, we adopt the Earth–Moon rotating frame described in Section \ref{subsec2.1} \cite{gong2007formation}. Differing from the origin defined as the Earth-Moon barycenter, the origin of the Earth-Moon rotating frame describing the relative motion is defined as the centroid of the spacecraft moving in the considered the NRHO (i.e., the chief spacecraft). The spacecraft perform formations around the chief spacecraft denotes the deputy spacecraft. The difference between the absolute states between the deputy and chief spacecraft denote the relative states $\delta \bm{X}=\left[\delta{x},\text{ }\delta{y},\text{ }\delta{z},\text{ }\delta{u},\text{ }\delta{v},\text{ }\delta{w}\right]^{\text{T}}$, which are defined as $\delta\bm{X}=\bm{X}_d-\bm{X}_{\text{NRHO}}$, where the subscript “\textit{d}” denotes the states associated with the deputy spacecraft. Let $\delta{x}$, $\delta{y}$, and $\delta{z}$ denote the relative coordinates, which are align with $x$-axis, $y$-axis, and $z$-axis defined in Section \ref{subsec2.1} respectively. Additionally, the local-vertical/local-horizontal (LVLH) frame can also be used for describing relative motion in the CR3BP, as discussed in Refs.  \cite{sanchez2020chance,chihang2023close}.

In the linear model, the initial relative states are defined as $\delta\bm{X}\left(t_0\right)=\bm{X}_d\left(t_0\right)-\bm{X}_{\text{NRHO}}\left(t_0\right)$, where the subscript “0” denotes the states associated with the initial epoch $t_0$. After propagating the orbits for a time interval $\Delta{t}$, the relative state is approximated by \cite{gong2007formation}:

\begin{equation}
\delta\bm{X}\left(t_0+\Delta{t}\right)\approx \bm{\Phi}\left(t_0+\Delta{t},\text{ }t_0\right)\delta\bm{X}\left(t_0\right)\label{eq55}
\end{equation}
where $\bm{\Phi}\left(t_0+\Delta{t},\text{ }t_0\right)$ denotes the state transition matrix (STM) from $t_0$ to $t_0+\Delta{t}$ and satisfies the following variational equation \cite{singh2021eclipse}:

\begin{equation}
\left\{ \begin{gathered}
  {{\bm{\dot \Phi }}} = \frac{{\partial \bm{f}}}{{\partial \bm{X}}}{\bm{\Phi }} \hfill \\
  {\bm{\Phi }}\left( t_0,\text{ }t_0 \right) = {\bm{I}_{6 \times 6}} \hfill \\ 
\end{gathered}  \right.\label{eq5}
\end{equation}
where $\bm{f}$ denotes the dynamical equations of the CR3BP Eq. \eqref{eq1}.

\subsection{Nonlinear Model for Relative Motion in the CR3BP}\label{subsec2.3}
While the linear model is useful for providing initial guesses, the strong nonlinearities in the CR3BP, especially near the perilune of the NRHO, necessitate the use of a nonlinear model for describing relative motion accurately \cite{jenson2024bounding,cunningham2023interpolated}. Therefore, this short communication develops methods for designing teardrop hovering formations along the NRHO in the context of the nonlinear model \cite{khoury2020orbital}. In this model, the relative states between the the deputy spacecraft and the chief spacecraft at any epoch is calculated by $\delta\bm{X}\left(t\right)=\bm{X}_d\left(t\right)-\bm{X}_{\text{NRHO}}\left(t\right)$, with both $\bm{X}_d\left(t\right)$ and $\bm{X}_{\text{NRHO}}\left(t\right)$ being propagated by Eq. \eqref{eq1}.

\section{Design and Continuation of Teardrop Hovering Formations along the NRHO}\label{sec3}
In this section, we develop methods for designing teardrop hovering formations along the NRHO. First, we introduce the concept of the teardrop hovering formation. Then, we present a method for designing teardrop hovering formations with relatively short revisit distances. Based on the results obtained from the first method, the method for designing hovering trajectories with relatively long revisit distances is developed.
\subsection{Teardrop Hovering Formation}\label{subsec3.1}
The teardrop hovering formation is a special type of hovering formation applicable to on-orbit service \cite{bai2020teardrop}. In the two-body dynamics, several methods for designing teardrop hovering formations have been proposed \cite{zhang2016flyaround,bai2020teardrop,bai2022practical}. In this short communication, the concept of the teardrop hovering formation is extended to the scenarios along the three-body periodic orbits. Differing from the two-body dynamics, the CR3BP lacks the analytical solution, which complicates the design of relative trajectories. Therefore, the purpose of this short communication is to develop numerical methods to design teardrop hovering formations along the NRHO and to provide further insights into future on-orbit service missions for targets moving in the NRHO.

The schematic of the teardrop hovering formation in the Earth-Moon rotating frame describing relative motion is shown in Fig. \ref{fig3}. In this formation, the deputy spacecraft revisits a fixed position relative to the chief spacecraft after a specific time interval $\Delta t$ (i.e., the revisit period). This fixed revisit position is denoted as $\bm{\rho}=\left[\delta{x}_r,\text{ }\delta{y}_r,\text{ }\delta{z}_r\right]$. To satisfy the revisit condition, an impulse $\Delta{\bm{v}}$ is performed at each revisit. To perform the teardrop hovering formation, the relative states should satisfy the following necessary condition \cite{bai2020teardrop,bai2022practical} ($\delta\bm{X}=\left[\delta\bm{r},\text{ }\delta\bm{v}\right]^{\text{T}}$):

\begin{figure}[h]
\centering
\includegraphics[width=0.3\textwidth]{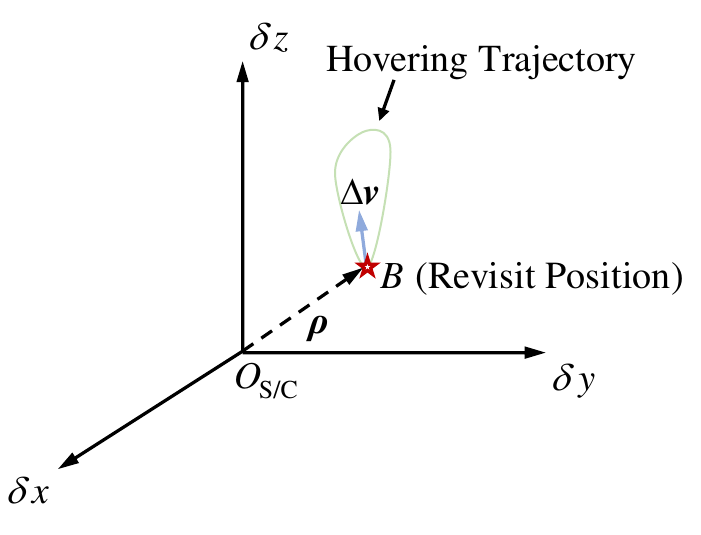}
\caption{The schematic of the teardrop hovering formation.}
\label{fig3}
\end{figure}

\begin{equation}
  \delta \bm{r}\left( {{t_j} + \Delta t} \right) - \delta \bm{r}\left( {{t_j}} \right) = \bm{0}, 
  {t_j} = {t_{j - 1}} + \Delta t,{\text{ }}j = 1,{\text{ }}2,{\text{ }}... 
\label{eq6}
\end{equation}
\begin{equation}
\delta\bm{r}\left(t_0\right)=\bm{\rho}\label{eq7}
\end{equation}

We considered a specific case of the teardrop hovering formation, namely, 1:1 hovering formation. In this scenario, $\Delta{t}=T_{\text{NRHO}}$, the constraints Eq. \eqref{eq7} can be specified as:

\begin{equation}\label{eq15}
\bm{\psi}  = \left[ \begin{gathered}
  \delta x\left( {{t_0} + {T_{{\text{NRHO}}}}} \right) - \delta x\left( {{t_0}} \right) \hfill \\
  \delta y\left( {{t_0} + {T_{{\text{NRHO}}}}} \right) - \delta y\left( {{t_0}} \right) \hfill \\
  \delta z\left( {{t_0} + {T_{{\text{NRHO}}}}} \right) - \delta z\left( {{t_0}} \right) \hfill \\ 
\end{gathered}  \right] = \bm{0}
\end{equation}
and the impulses performed for each revisit are identical (denoted as $\Delta{\bm{v}}$ and $\Delta{v}$) and can be calculated by:

\begin{equation}
\Delta{\bm{v}}=\delta{\bm{v}}\left(t_0\right)-\delta{\bm{v}}\left(t_0+T_{\text{NRHO}}\right)\label{eq888}
\end{equation}
\begin{equation}
\Delta{v}=\left|\left|\delta{\bm{v}}\left(t_0\right)-\delta{\bm{v}}\left(t_0+T_{\text{NRHO}}\right)\right|\right|\label{eq8}
\end{equation}

In the following text, we detail the methods for designing and continuing this type of teardrop hovering formation along the considered NRHO.

\subsection{Design of Teardrop Hovering Formations with Relatively Short Revisit Distances}\label{subsec3.2}
\begin{figure}[h]
\centering
\includegraphics[width=0.25\textwidth]{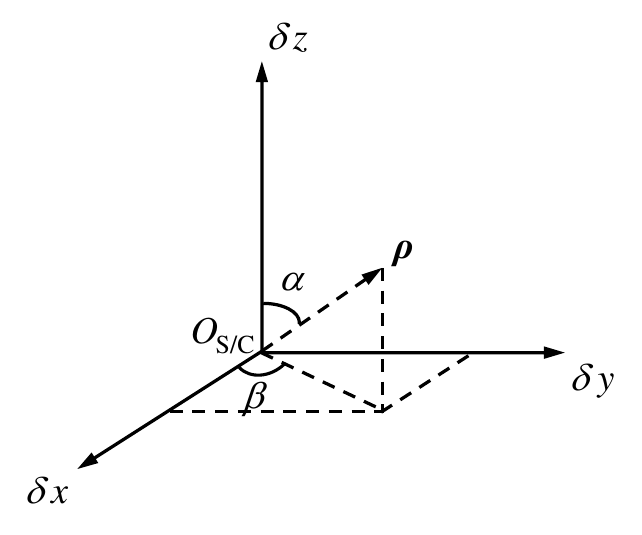}
\caption{The geometry of $\alpha$ and $\beta$.}
\label{fig4}
\end{figure}

We first focus on designing the teardrop hovering formations with relatively short revisit distances. The initial relative position is set equal to the revisit position, i.e., $\delta\bm r\left(t_0\right)=\bm{\rho}$. We adopt the fixed value of $\rho$ ($\rho=\left|\left|\bm{\rho}\right|\right|$ denotes the revisit distance) in this case, i.e., $\rho=1\text{ km}$. Therefore, the components of $\bm\rho$ can be expressed as:
\begin{equation}
\left\{ \begin{gathered}
  \delta {x_r} = \rho \sin \alpha \cos \beta  \hfill \\
  \delta {y_r} = \rho \sin \alpha \sin \beta  \hfill \\
  \delta {z_r} = \rho \cos \alpha  \hfill \\ 
\end{gathered}  \right.\label{eq9}
\end{equation}
where the geometry of $\alpha$ and $\beta$ can be found in Fig. \ref{fig4}.

Since the revisit distance is short, the nonlinear effects on relative motion are relatively minor. Therefore, we use the linear model to generate initial guesses of the hovering trajectories, which are then corrected to satisfy the revisit conditions in the full nonlinear model. Given the initial relative position, the design variable is the initial relative velocity $\delta{\bm{v}}\left(t_0\right)$. In the linear model, the relative state after $T_{\text{NRHO}}$ can be calculate by:
\begin{equation}
\left[ {\begin{array}{*{20}{c}}
  {\delta \bm{r}\left( {{t_0} + {T_{{\text{NRHO}}}}} \right)} \\ 
  {\delta \bm{v}\left( {{t_0} + {T_{{\text{NRHO}}}}} \right)} 
\end{array}} \right] \approx \left[ {\begin{array}{*{20}{c}}
  {{\bm{\Phi }_{rr}}}&{{\bm{\Phi }_{rv}}} \\ 
  {{\bm{\Phi }_{vr}}}&{{\bm{\Phi }_{vv}}} 
\end{array}} \right]\left[ {\begin{array}{*{20}{c}}
  {\delta \bm{r}\left( {{t_0}} \right)} \\ 
  {\delta \bm{v}\left( {{t_0}} \right)} 
\end{array}} \right]
\label{eq10}
\end{equation}
\begin{equation}\label{eq11}
\bm{\Phi}=\bm{\Phi}\left(t_0+T_{\text{NRHO}},\text{ }t_0\right)=\left[ {\begin{array}{*{20}{c}}
  {{\bm{\Phi }_{rr}}}&{{\bm{\Phi }_{rv}}} \\ 
  {{\bm{\Phi }_{vr}}}&{{\bm{\Phi }_{vv}}} 
\end{array}} \right]\
\end{equation}
The STM shown in Eq. \eqref{eq11} denotes the STM associated with the considered NRHO (i.e., the monodromy matrix of the NRHO). Therefore, the relative position at $t_0+T_{\text{NRHO}}$ can be expressed as \cite{gong2007formation}:
\begin{equation}\label{eq12}
\delta \bm{r}\left( {{t_0} + {T_{{\text{NRHO}}}}} \right) \approx {\bm {\Phi }_{rr}}\delta \bm {r}\left( {{t_0}} \right) + {\bm {\Phi }_{rv}}\delta \bm {v}\left( {{t_0}} \right)
\end{equation}
To satisfy the revisit conditions, the initial relative velocity should be chosen such that the relative position at $t_0+T_{\text{NRHO}}$ matches the initial position. The initial guess of $\delta \bm {v}\left( {{t_0}} \right)$ is calculated by:
\begin{equation}\label{eq13}
{\bm {\Phi }_{rv}}\delta \bm {v}\left( {{t_0}} \right) \approx\delta \bm{r}\left( {{t_0} + {T_{{\text{NRHO}}}}} \right)-{\bm {\Phi }_{rr}}\delta \bm {r}\left( {{t_0}} \right) =\left(\bm{I}_{3\times 3}-{\bm {\Phi }_{rr}}\right)\delta \bm {r}\left( {{t_0}} \right)
\end{equation}
\begin{equation}\label{eq14}
\delta \bm{v}\left( {t_0}  \right) = {\bm {\Phi }_{rv}}^{-1}\left(\bm{I}_{3\times 3}-{\bm {\Phi }_{rr}}\right)\delta \bm {r}\left( {{t_0}} \right)
\end{equation}
In the simulations, the inverse of the matrix ${\bm {\Phi }_{rv}}$ is achieve by MATLAB®'s pinv command. Once the initial guess of $\delta \bm{v}\left( {t_0}  \right)$ is determined, we correct it to the nonlinear model. The correction is transformed into a nonlinear programming (NLP) problem, performed by MATLAB®'s fmincon command. In the NLP problem, the optimization variable is $\delta \bm{v}\left( {t_0}  \right)$, and the constraints of the NLP problem can be specifically expressed as Eq. \eqref{eq15} and the function to be minimized is set as $\left\| \bm{\psi} \right\|$. The parameter settings for the fmincon command used in simulations are provided in Appendix \ref{secA1}. After optimization, the impulse is calculated by Eqs. \eqref{eq888}-\eqref{eq8}.

\textit{Remark 1}: In this short communication, we focus on the 1:1 teardrop hovering formation. In this case, due to the periodicity of the NRHO ($\bm{X}_{\text{NRHO}}\left(t_0+T_{\text{NRHO}}\right)=\bm{X}_{\text{NRHO}}\left(t_0\right)$), $\delta \bm r\left(t_0+T_{\text{NRHO}}\right)-\delta \bm r\left(t_0\right)$ is equivalent to $\bm{r}\left(t_0+T_{\text{NRHO}}\right)-\bm{r}\left(t_0\right)$, where $\bm{X}=\left[\bm{r}^{\text{T}},\text{ }\bm{v}^{\text{T}}\right]^{\text{T}}$ denotes the states of the absolute hovering trajectories. Meanwhile, the impulse vector can be calculated by $\bm{v}\left(t_0\right)-\bm{v}\left(t_0+T_{\text{NRHO}}\right)$.

\textit{Remark 2}: Compared to the teardrop hovering formations designed solely using the linear model, the formations corrected by our method are expected to achieve significantly lower errors in the revisit position after each revisit period, as discussed in detail in Section \ref{subsec4.1}.

\textit{Remark 3}: Using the parameterization shown in Eq. \eqref{eq9} and the method developed in this subsection, a teardrop hovering formation can be designed from grid-sampled values of $\left(\alpha,\text{ }\beta\right)$. Moreover, the resulting impulse distribution with respect to $\left(\alpha,\text{ }\beta\right)$ can provide a valuable reference for selecting the optimal values of $\left(\alpha,\text{ }\beta\right)$ and revisit position in practical on-orbit service mission designs.

Subsequently, based on the design results of teardrop hovering formations with relatively short revisit distances, we further propose the method for designing and continuing teardrop hovering formations with relatively long revisit distances.

\subsection{Design and Continuation of Teardrop Hovering Formations with Relatively Long Revisit Distances}\label{subsec3.3}
In this subsection, we address the design of teardrop hovering formations with relatively long revisit distances. When the revisit distance is increased, nonlinear effects become more noticeable, and the method from Section \ref{subsec3.2} may no longer yield satisfactory performance. To generate more accurate initial guesses in this case, we develop a new form of numerical continuation method. However, for this new scenario, the continuation scenario can be different from those reported in previous works \cite{oshima2019low,oshima2022continuation,singh2019mission,alvarado2024exploration}. The continuation scenario requiring to address in this short communication belongs to trajectory continuation in a parameter that is part of the constraints, i.e., the parameter $\rho$. This scenario has not been thoroughly investigated, and a new continuation method has been developed to meet the new design requirements.

Typically, the continuation method is a predictor-corrector method. In this case, specifically, the feasible direction of the constrained trajectories (the teardrop hovering formation in this short communication) is predicted using the existing trajectories obtained in Section \ref{subsec3.2}. This prediction is then corrected to obtain the next feasible solution. In this short communication, we adopt the linear predictor to continue hovering trajectories in the revisit distance $\rho$. Here, the initial revisit position from Section \ref{subsec3.2} is denoted as $\bm{\rho}^0$. As the value of $\rho$ varies, the new positions are denoted as $\bm\rho^{k}$ for $\text{ }k=1,\text{ }2,\text{ }...$, where $\bm\rho^{k}=\bm\rho^{k-1}+\Delta{\bm\rho}$ with the step-size $\Delta\rho$ set to 0.1 km. The components of variation $\Delta{\bm\rho}$ can be expressed as:
\begin{equation}
\left\{ \begin{gathered}
  \Delta\left(\delta {x_r}\right) = \Delta\rho \sin \alpha \cos \beta  \hfill \\
  \Delta\left(\delta {y_r}\right) = \Delta\rho \sin \alpha \sin \beta  \hfill \\
  \Delta\left(\delta {z_r}\right) = \Delta\rho \cos \alpha  \hfill \\ 
\end{gathered}  \right.\label{eq17}
\end{equation}
As an example, consider the first predictor–corrector step from $\bm{\rho}^0$ to $\bm{\rho}^1$. The difference in constraints Eq. \eqref{eq15} can be approximated by the first-order Taylor expansion \cite{fu2025analytical}:
\begin{align}\label{eq18}
{\bm{\psi }^1} - {\bm{\psi }^0} &\approx \left.\frac{{\partial \bm{\psi }}}{{\partial \left( {\delta \bm{r}\left( {{t_0}} \right)} \right)}}\right|_{{\bm{\rho } = {\bm{\rho }^0}}}\Delta \left( {\delta \bm{r}\left( {{t_0}} \right)} \right) 
+ \left.\frac{{\partial \bm{\psi }}}{{\partial \left( {\delta \bm{v}\left( {{t_0}} \right)} \right)}}\right|_{{\bm{\rho } = {\bm{\rho }^0}}}\Delta \left( {\delta \bm{v}\left( {{t_0}} \right)} \right) \\ 
  \notag & = \left.\frac{{\partial \bm{\psi }}}{{\partial \left( {\delta \bm{r}\left( {{t_0}} \right)} \right)}}\right|_{{\bm{\rho } = {\bm{\rho }^0}}}\Delta \bm{\rho } 
 + \left.\frac{{\partial \bm{\psi }}}{{\partial \left( {\delta \bm{v}\left( {{t_0}} \right)} \right)}}\right|_{{\bm{\rho } = {\bm{\rho }^0}}}\Delta \left( {\delta \bm{v}\left( {{t_0}} \right)} \right) \\ 
  \notag & = \bm{0} 
\end{align}
Equation \eqref{eq18} can be reconstructed in the following matrix form:
\begin{equation}
\bm{A}\Delta \left( {\delta \bm{v}\left( {{t_0}} \right)} \right) = \bm{b}\label{eq19}
\end{equation}
For the 1:1 hovering formation along the considered NRHO, the aforementioned equations can be further simplified. Since the teardrop hovering formation is performed along the same NRHO with the same revisit period (i.e., $\Delta t=T_{\text{NRHO}}$), the differences in relative states at $t_0$ and $t_0+T_{\text{NRHO}}$ are equivalent to the differences in absolute states of the hovering trajectory (the absolute states of the hovering trajectory are expressed as $\bm{X}_{\text{NRHO}}+\delta\bm{X}$). Therefore, the aforementioned partial derivatives can be further specified as the partial derivatives associated with the absolute states. The components of the matrix $\bm{A}$ and the vector $\bm{b}$ can be expressed as:
\begin{equation}
\begin{gathered}
  {\bm{A}_{11}} = {\bm{\Phi }_{14}}{\text{ }\text{ }}{\bm{A}_{12}} = {\bm{\Phi }_{15}}{\text{ }\text{ }}{\bm{A}_{13}} = {\bm{\Phi }_{16}} \hfill \\
  {\bm{A}_{21}} = {\bm{\Phi }_{24}}{\text{ }\text{ }}{\bm{A}_{22}} = {\bm{\Phi }_{25}}{\text{ }\text{ }}{\bm{A}_{23}} = {\bm{\Phi }_{26}} \hfill \\
  {\bm{A}_{31}} = {\bm{\Phi }_{34}}{\text{ }\text{ }}{\bm{A}_{32}} = {\bm{\Phi }_{35}}{\text{ }\text{ }}{\bm{A}_{33}} = {\bm{\Phi }_{36}} \hfill \\ 
\end{gathered}
\label{eq26}
\end{equation}

\begin{equation}\label{eq27}
\begin{gathered}
  {\bm{b}_1} = \left( {1 - {\bm{\Phi }_{11}}} \right)\Delta \left( {\delta {x_r}\left( {{t_0}} \right)} \right)  - {\bm{\Phi }_{12}}\Delta \left( {\delta {y_r}\left( {{t_0}} \right)} \right)  - {\bm{\Phi }_{13}}\Delta \left( {\delta {z_r}\left( {{t_0}} \right)} \right) \hfill \\
  {\bm{b}_2} = \left( {1 - {\bm{\Phi }_{22}}} \right)\Delta \left( {\delta {y_r}\left( {{t_0}} \right)} \right)  - {\bm{\Phi }_{21}}\Delta \left( {\delta {x_r}\left( {{t_0}} \right)} \right)  - {\bm{\Phi }_{23}}\Delta \left( {\delta {z_r}\left( {{t_0}} \right)} \right) \hfill \\
  {\bm{b}_3} = \left( {1 - {\bm{\Phi }_{33}}} \right)\Delta \left( {\delta {z_r}\left( {{t_0}} \right)} \right)  - {\bm{\Phi }_{31}}\Delta \left( {\delta {x_r}\left( {{t_0}} \right)} \right)  - {\bm{\Phi }_{32}}\Delta \left( {\delta {y_r}\left( {{t_0}} \right)} \right) 
\end{gathered}
\end{equation}
where the STM $\bm{\Phi}$ denotes the STM associated with the absolute hovering trajectory along the NRHO from $t_0$ to $t_0+T_{\text{NRHO}}$, which is different from the STM mentioned in Section \ref{subsec3.2}. Given the components of the matrix $\bm{A}$ and the vector $\bm{b}$, the prediction of $\Delta \left( {\delta \bm{v}\left( {{t_0}} \right)} \right)$ is generated by the least squares method \cite{wu2022analytical}:
\begin{equation}
\Delta \left( {\delta \bm{v}\left( {{t_0}} \right)} \right)=\left(\bm{A}^{\text{T}}\bm{A}\right)^{-1}\bm{A}^{\text{T}}\bm{b}\label{eq30}
\end{equation}
This solution provides a prediction for the required adjustment $\Delta \left( {\delta \bm{v}\left( {{t_0}} \right)} \right)$ in the initial relative velocity. This updated velocity is used as the initial guess for the predictor–corrector step, and subsequent differential corrections ensure that the constraints Eq. \eqref{eq15} are satisfied. The process is repeated until the termination criteria are satisfied as follows:

\begin{enumerate}[label=(\arabic*)]
\item The modulus of the constraint equations $\bm{\psi}=\textbf{0}$ after differential correction exceeds the tolerance $\varepsilon_{\text{tol}}$ set to $1\times 10^{-9}$, i.e., $\left\| \bm{\psi} \right\| > 1\times 10^{-9}$;
\item The number of predictor-corrector steps ($j$) exceeds the maximum number of predictor-corrector steps $N$ set to 499, i.e., $j>499$.
\end{enumerate}

\begin{figure}[h]
\centering
\includegraphics[width=0.4\textwidth]{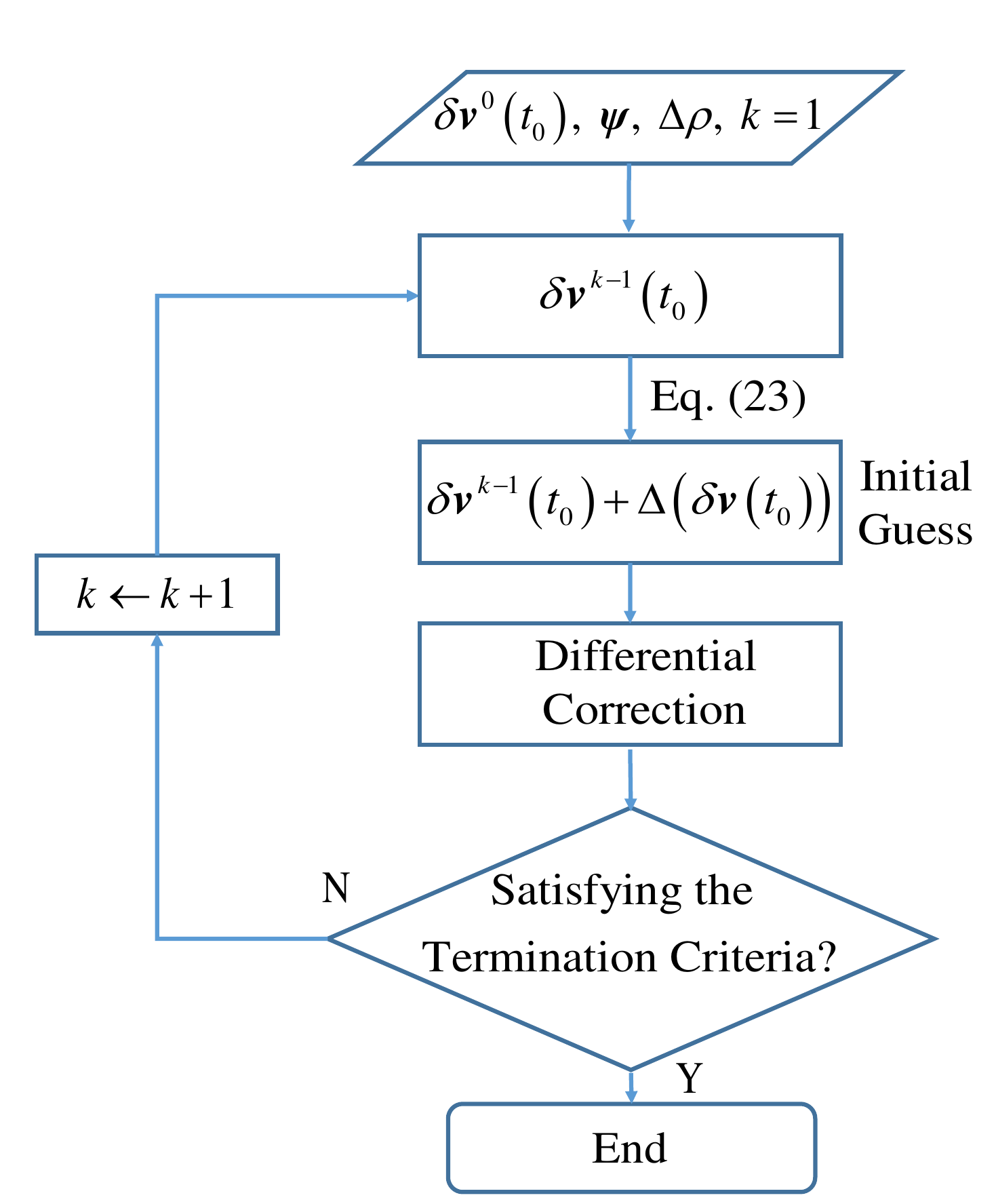}
\caption{The flow chart of the continuation method.}
\label{fig5}
\end{figure}

Then, the procedure of the continuation method developed in this short communication is presented in Fig. \ref{fig5}.

\textit{Remark 4}: For the new continuation scenario, i.e., continuation in a parameter that is a part of the constraints, a new form of linear predictor has been developed to meet the design requirements. As the value of $\rho$ increases, the nonlinear effects become more significant, leading to larger differences between the linear and nonlinear models in the design results.

\textit{Remark 5}: This short communication focuses on designing teardrop hovering formations along the NRHO. However, the methods developed in Sections \ref{subsec3.2}-\ref{subsec3.3} can be extended to design teardrop hovering formations along other multi-body periodic orbits (e.g., DROs).

\textit{Remark 6}: This short communication presents the design of a teardrop hovering formation in the Earth–Moon rotating frame and clarifies that the same design methods can be used in the LVLH frame.

In the next section, we present the design results and discuss the necessity of using the nonlinear model for designing teardrop hovering formations along the considered NRHO.

\section{Results and Discussion}\label{sec4}
In this section, the design results obtained from the methods developed in Section \ref{subsec3.2}-\ref{subsec3.3} are presented. The comparison between results obtained from the methods using the linear and nonlinear models is then performed. Here the method using the linear models refers to $\delta \bm{v}\left(t_0\right)$ and $\Delta \bm{v}$ generated from the linear model without differential correction.
\subsection{Teardrop Hovering Formations with Relatively Short Revisit Distances}\label{subsec4.1}
\begin{figure}[h]
\centering
\includegraphics[width=0.3\textwidth]{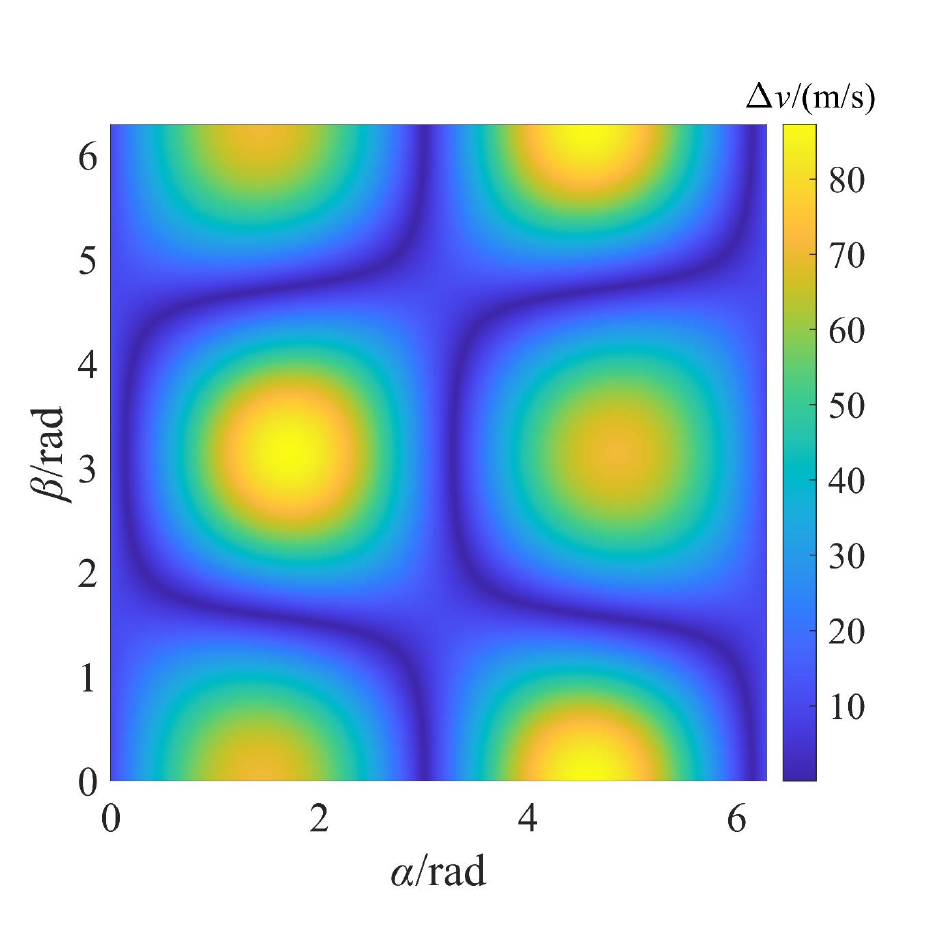}
\caption{The distributions of the impulse for each pair $\left(\alpha,\text{ }\beta\right)$.}
\label{fig6}
\end{figure}

For this case, we set $\alpha$ and $\beta$ as $\alpha\in\left[0,\text{ }2\pi\right]$ and $\beta\in\left[0,\text{ }2\pi\right]$ with a step-size of $\pi/100$, and the value of $\rho$ is fixed at 1 km. At $t=t_0$, the states of the chief spacecraft (i.e., the states of the NRHO) are taken as given in Table \ref{tab2} (other initial states of the the chief spacecraft moving in the NRHO can be selected for specific mission requirement). For each pair $\left(\alpha,\text{ }\beta\right)$, the corresponding teardrop hovering formation along the considered NRHO is designed by the method developed in Section \ref{subsec3.2}. The impulse distribution in this case is presented in Fig. \ref{fig6}. This impulse distribution can serve as a valuable reference for selecting the values of $\left(\alpha,\text{ }\beta\right)$ for the design of the teardrop hovering formation.

\begin{figure}[h]
\centering
\includegraphics[width=0.9\textwidth]{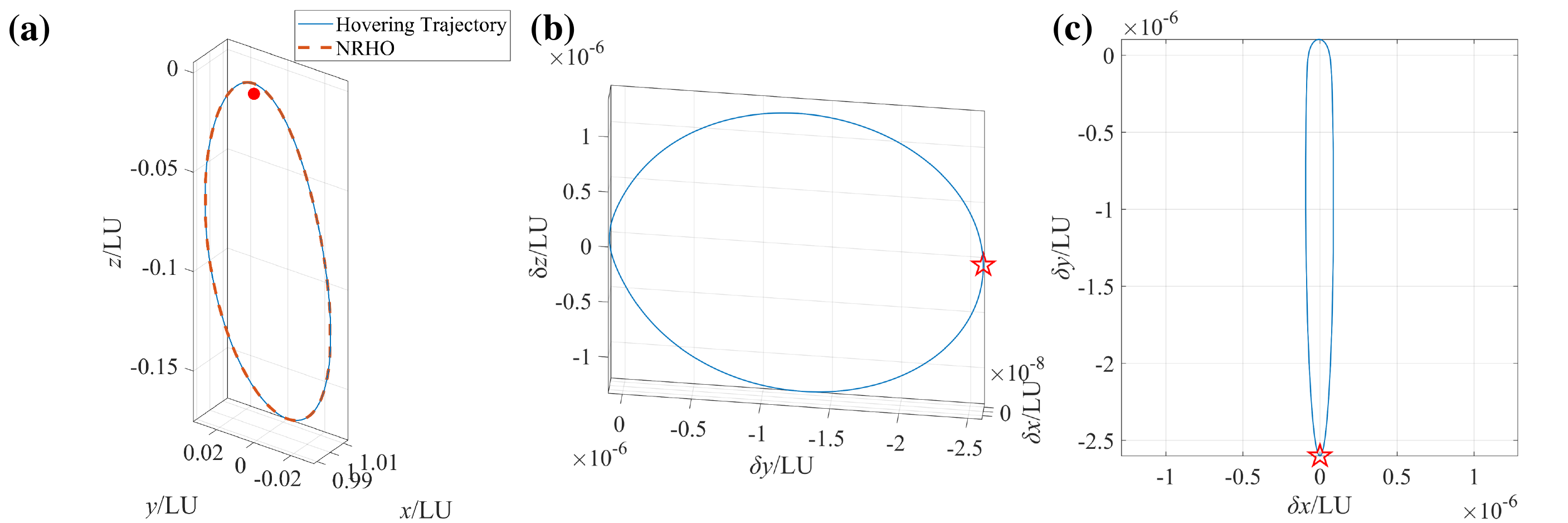}
\caption{The absolute trajectory and relative trajectory of the selected example during one revisit period. (a) Absolute trajectory; (b) Relative trajectory; (c) Relative trajectory in the $\delta x-\delta y$ plane. The red pentagram denotes the revisit position.}
\label{fig7}
\end{figure}

As shown in Fig. \ref{fig6}, the impulse distribution exhibits symmetry with respect to $\beta=\pi$. Meanwhile, the maximum value of impulse to satisfy the revisit conditions is relatively high even when $\rho=1\text{ km}$. To ensure engineering practicality, we select the trajectory that requires the minimum impulse for further analysis. Notably, this trajectory with minimum impulse is obtained from the numerical search. The relative states of this trajectory at $t=t_0$ are presented in Table \ref{tab3}. The minimum impulse required is $7.333\times 10^{-4}$ m/s. Given that the revisit period is approximately 6-7 days, we consider this trajectory suitable for practical on-orbit service missions. Figure \ref{fig7} presents both the absolute and relative trajectories of this example during one revisit period. 

\begin{table}
\caption{The relative states of the hovering trajectory with the minimum impulse at $t=t_0$}\label{tab3}%
\centering
\renewcommand{\arraystretch}{1.5}
\begin{tabular}{@{}lll@{}}
\hline
Parameter & Value  & Unit \\
\hline
$\delta x\left(t_0\right)$ & $0$ & LU \\
$\delta y\left(t_0\right)$ & $-2.60142297836917 \times 10^{-6}$ & LU \\
$\delta z\left(t_0\right)$ & $0$ & LU \\
$\delta u\left(t_0\right)$ & $-3.2643727501816\times 10^{-5}$ & LU/TU \\
$\delta v\left(t_0\right)$ & $-1.98390221419\times 10^{-7}$ & LU/TU \\
$\delta w\left(t_0\right)$ & $5.33425501523417\times 10^{-4}$ & LU/TU \\
\hline
\end{tabular}
\end{table}

\begin{figure}[!htb]
\centering
\includegraphics[width=0.6\textwidth]{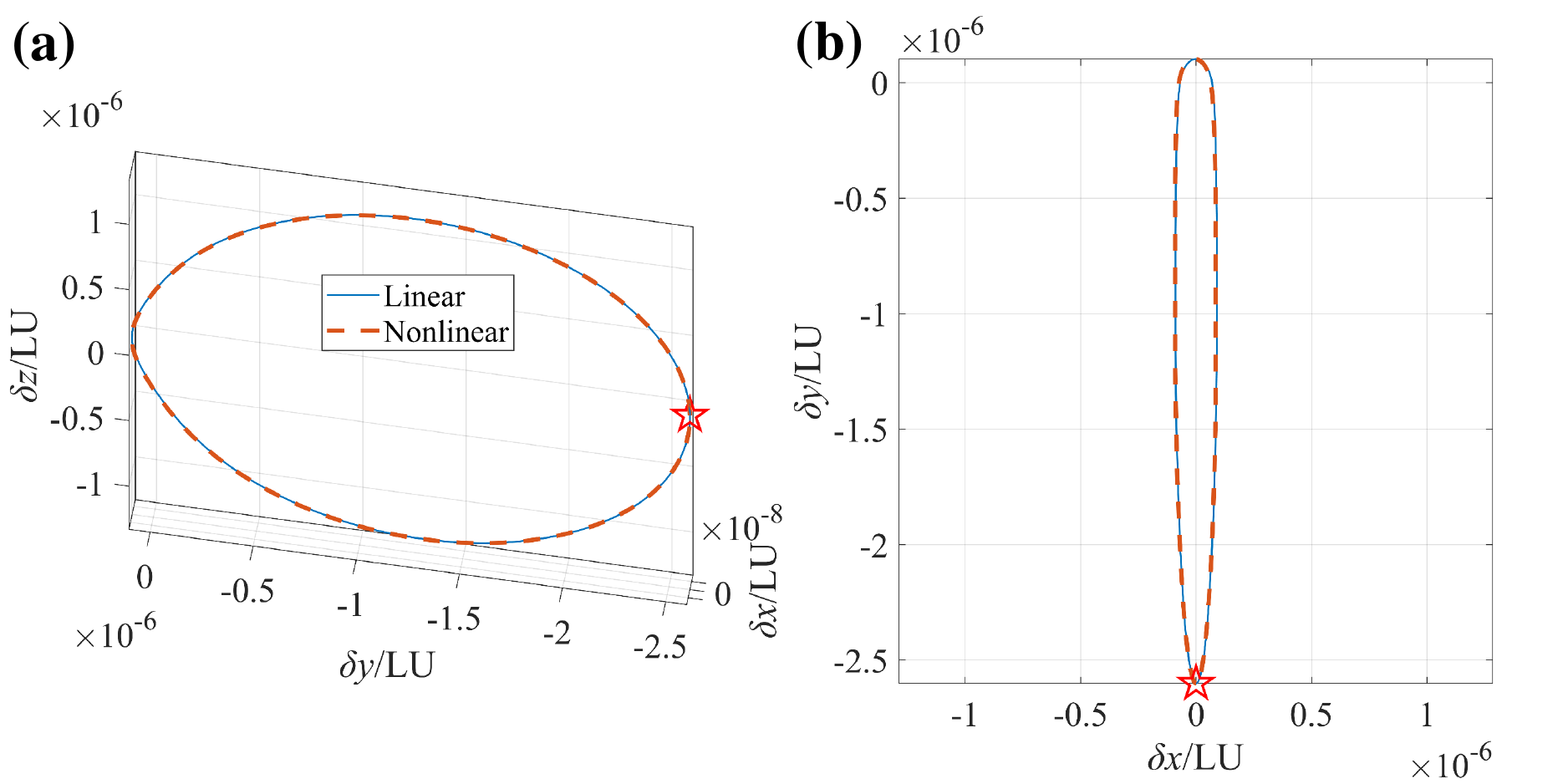}
\caption{The relative trajectories with the minimum impulse before and after differential correction. (a) Relative trajectories; (b) Relative trajectories in the $\delta x-\delta y$ plane.}
\label{fig13}
\end{figure}

In particular, for the $\left(\alpha,\text{ }\beta\right)$ pair that yields the minimum impulse, the vector $\delta{\bm{X}}_0$ combining $\delta{\bm{r}}_0$ with the $\delta{\bm{v}}_0$ generated from Eq. \eqref{eq14} is nearly aligned with the eigenvector of the NRHO’s monodromy matrix corresponding to the eigenvalue 1, i.e., the impulse calculated by the linear model is almost 0. Therefore, the actual impulse, $\Delta{v}=7.333\times 10^{-4}\text{ m/s}$, is obtained from considering the nonlinear model. The corresponding relative trajectories before and after differential correction (i.e., the teardrop hovering formation design by the linear and nonlinear models) is presented in Fig. \ref{fig13}. Linking the initial states presented in Table \ref{tab3}, we note that the absolute trajectory of the deputy spacecraft deviates from the NRHO, indicating that the teardrop hovering formation with the minimum impulse is a near-natural 1:1 formation achieved by exploiting the dynamical properties near the NRHO (not the natural teardrop hovering formation). The following analysis focus on this example.  

\begin{figure}[!htb]
\centering
\includegraphics[width=0.6\textwidth]{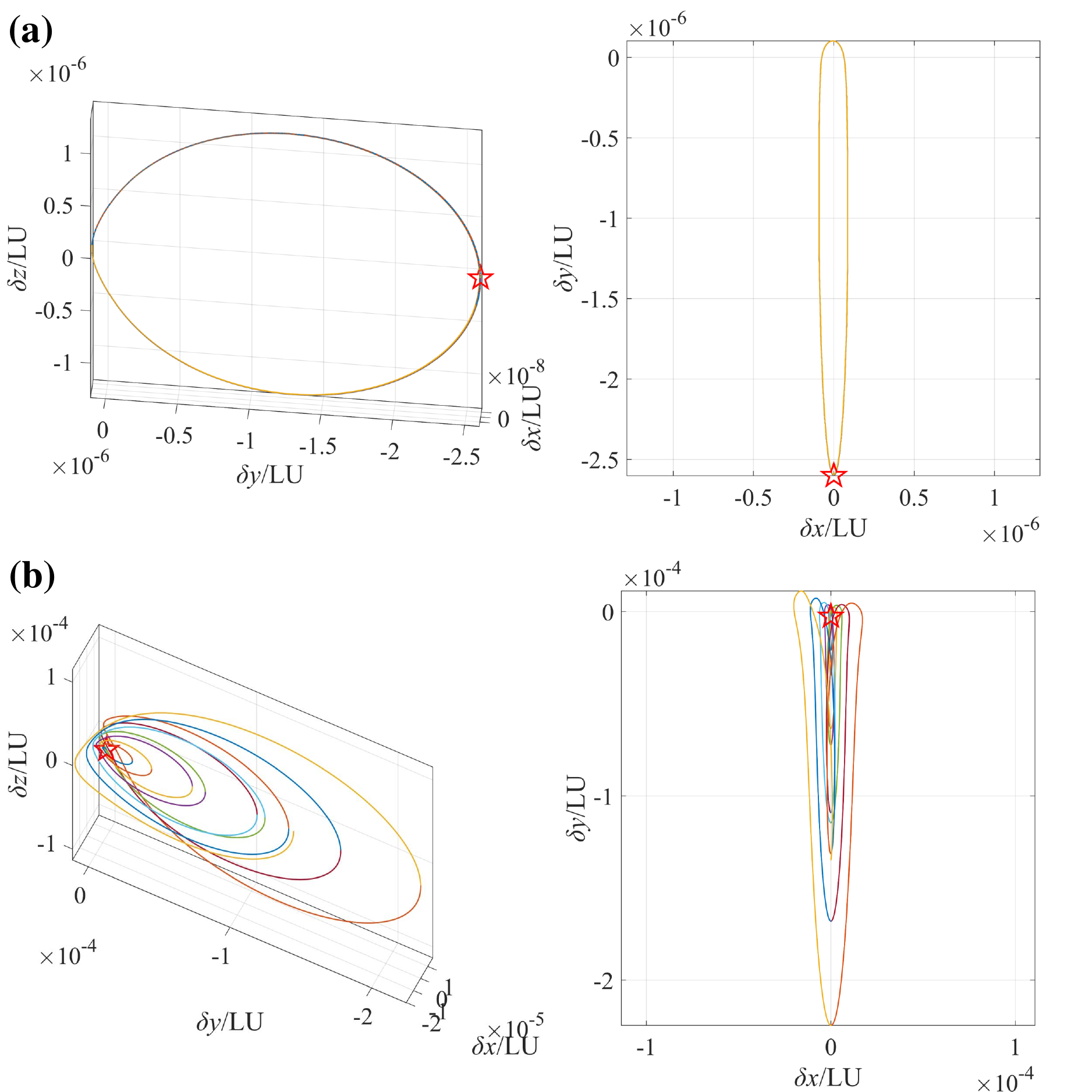}
\caption{The actual relative trajectories designed by linear and nonlinear models during 10 revisit periods. (a) The actual relative trajectories designed by nonlinear model; (b) The actual relative trajectories designed by linear model.}
\label{fig8}
\end{figure}

Although the configurations of relative trajectories shown in Fig. \ref{fig13} appear similar, the design results ($\Delta{\bm{v}}$ and $\delta{\bm{v}}_0$) obtained from the linear and nonlinear models can differ significantly over long-term hovering. Therefore, to further demonstrate the necessity of using the nonlinear model to design teardrop hovering formations along the considered NRHO, we compare the results obtained from the linear and nonlinear models. Using the same revisit position as specified in Table \ref{tab3}, we propagate the absolute states generated from the relative states obtained from both models (i.e., $\delta \bm{v}\left(t_0\right)$) using Eq. \eqref{eq1}. We then subtract the states of NRHOs to obtain the actual relative states. Impulses are applied at the end of each revisit period, and the impulses are calculated in the linear and nonlinear models, respectively. The actual relative trajectories designed by linear and nonlinear models during 10 revisit periods are presented in Fig. \ref{fig8}.

\begin{figure*}[!htb]
\centering
\includegraphics[width=0.6\textwidth]{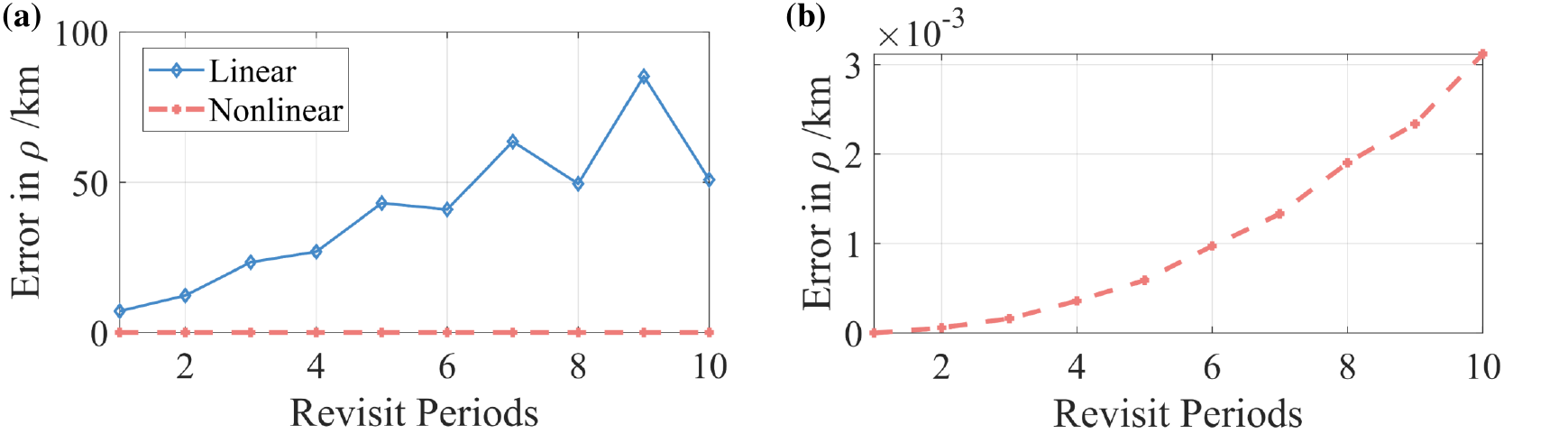}
\caption{The drifts relative to the given revisit position after each revisit period by linear and nonlinear models. (a) The drifts relative to the given revisit position after each revisit period by linear and nonlinear models; (b) The drifts relative to the given revisit position after each revisit period by the nonlinear model.}
\label{fig9}
\end{figure*}

From Fig. \ref{fig8}, it is observed that when examining the actual relative motion, the teardrop hovering formation designed by the linear model exhibits noticeable drift relative to the revisit position, even with the relatively low value of $\rho$. In contrast, the formation designed using the nonlinear model shows only a slight drift. This difference arises because the numerical methods used to compute $\delta \bm{v}\left(t_0\right)$ and $\Delta \bm{v}$ cannot enforce the constraint $\bm{\psi}=\bm{0}$ exactly. Figure \ref{fig9} quantitatively records the drifts (i.e., the errors in $\rho$) after each revisit period. These results indicate that the linear model introduces significant errors in the actual relative motion, whereas the nonlinear model yields much smaller errors, strengthening the necessity of using the nonlinear model.

\subsection{Teardrop Hovering Formations with Relatively Long Revisit Distances}\label{subsec4.2}
\begin{figure}[!htb]
\centering
\includegraphics[width=0.6\textwidth]{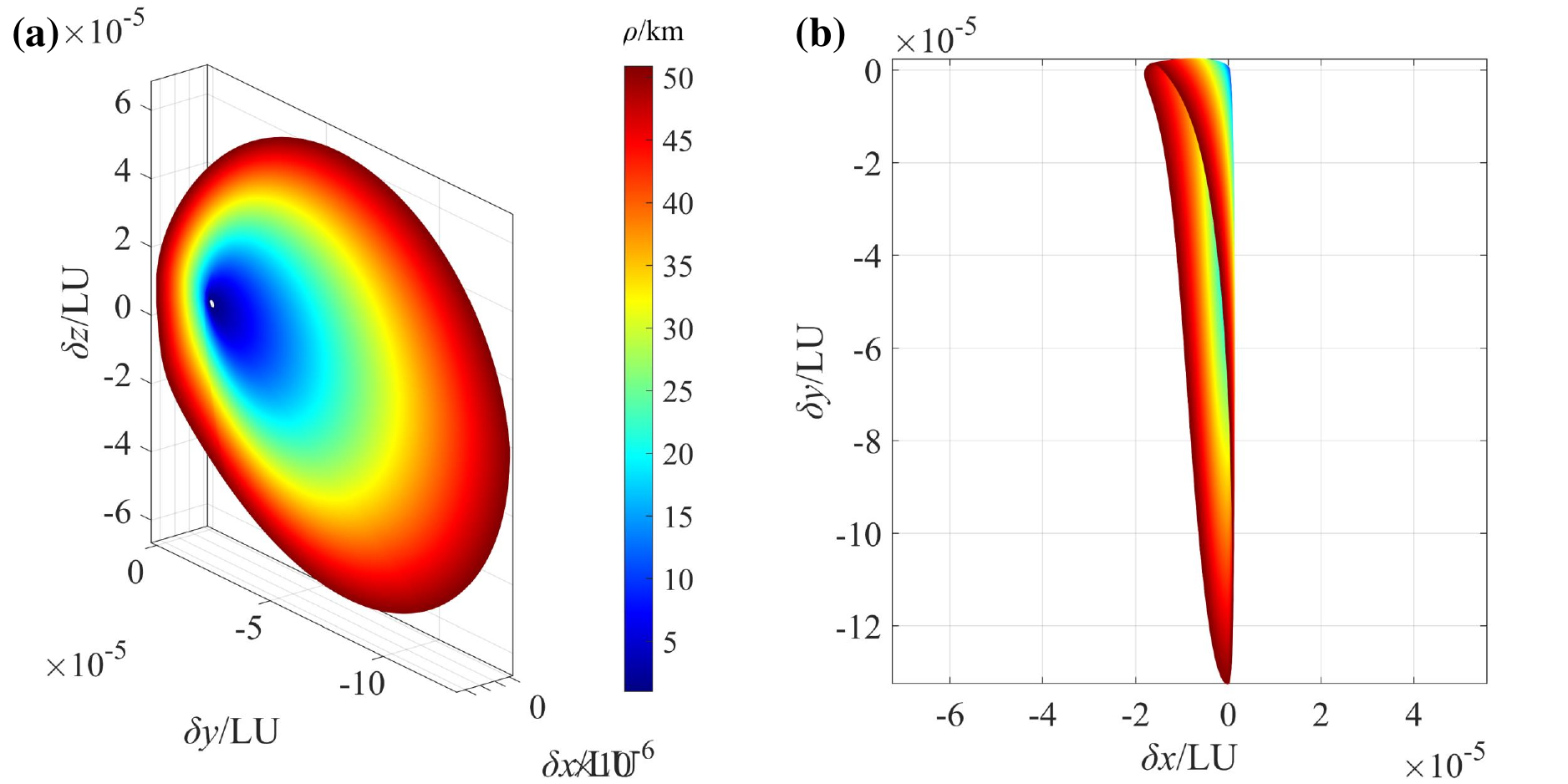}
\caption{The continued trajectories. (a) The continued trajectories; (b) The continued trajectories in the $\delta x-\delta y$ plane.}
\label{fig10}
\end{figure}

\begin{figure}[!htb]
\centering
\includegraphics[width=0.6\textwidth]{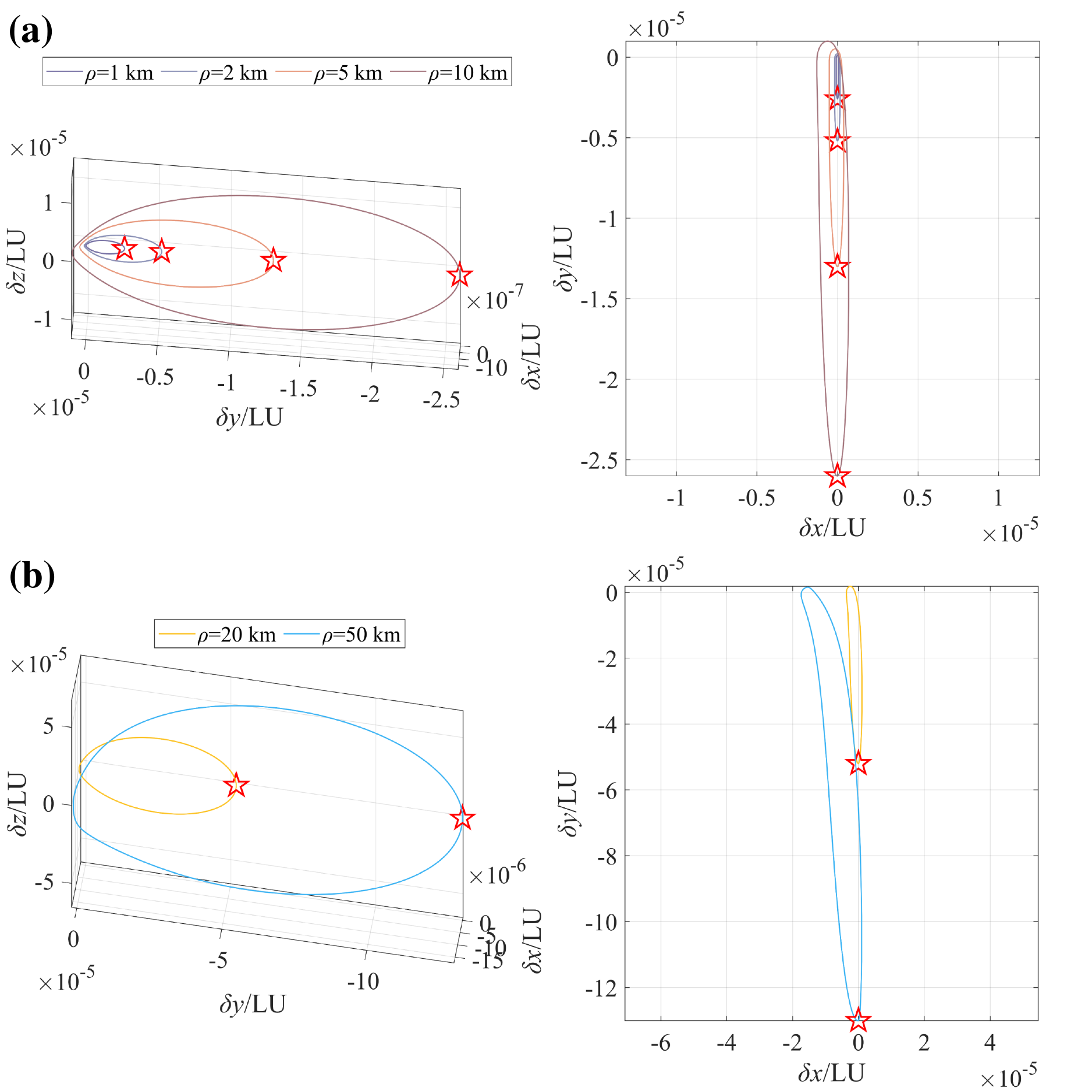}
\caption{Samples of the continued trajectories. (a) The continued trajectories with $\rho=1,\text{ }2,\text{ }5,\text{ }10\text{ km}$; (b) The continued trajectories with $\rho=20,\text{ }50\text{ km}$. The red pentagrams denote the different revisit positions.}
\label{fig11}
\end{figure}

\begin{figure}[!htb]
\centering
\includegraphics[width=0.6\textwidth]{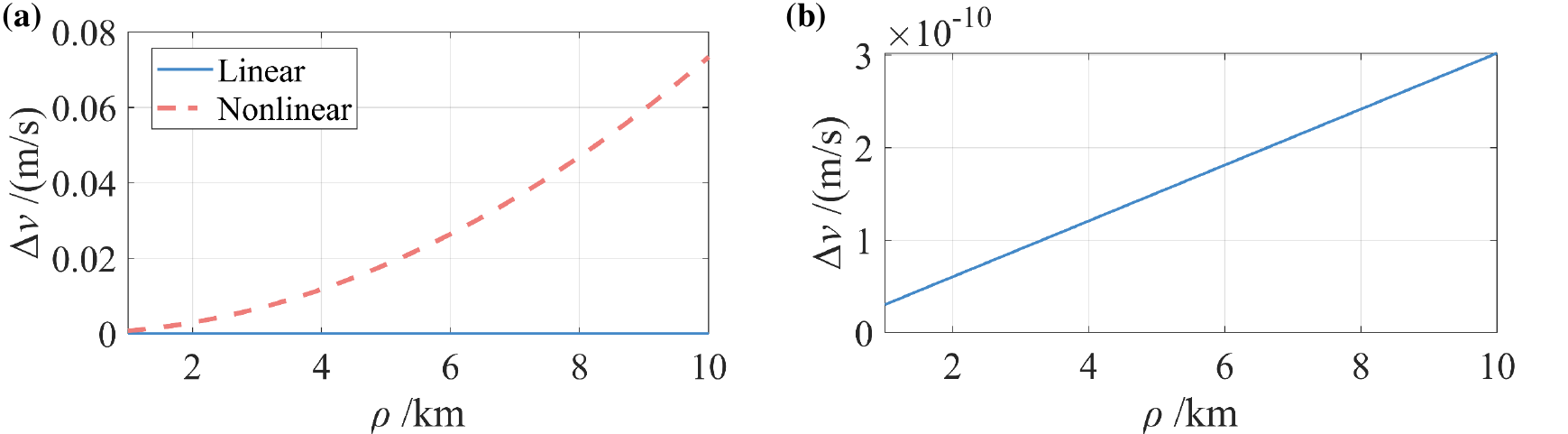}
\caption{The distributions of solutions expressed in terms of designed $\Delta v$ and $\rho$. (a) The distributions of solutions obtained from the linear and nonlinear models; (b) The distributions of solutions obtained from the nonlinear models.}
\label{fig12}
\end{figure}

In this subsection, we select the trajectory sample presented in Section \ref{subsec4.1} and continue it in the revisit distance from $\rho=1\text{ km}$ to higher values, with $\Delta \rho=0.1 \text{ km}$ and the same $\left(\alpha,\text{ }\beta\right)$. The continued trajectories with different values of $\rho$ are presented in Fig. \ref{fig10}. Some examples obtained from the developed continuation method is presented in Fig. \ref{fig11}. The distributions of solutions expressed in terms of $\Delta v$ and $\rho$ are presented in Fig. \ref{fig12}. Differing from the linear variation in $\Delta v$ predicted by the linear model, the impulses calculated using the nonlinear model are relatively high and exhibit a nonlinear relationship with respect to $\rho$. These results further demonstrate the errors introduced by the linear method and strengthen the necessity of using the nonlinear model.

\section{Conclusion}\label{sec5}
This short communication is devoted to the design of teardrop hovering formations along a Near Rectilinear Halo Orbit (NRHO). The 9:2 NRHO in the Earth-Moon circular restricted three-body problem (CR3BP) is selected as the reference orbit. We extend the concept of the teardrop hovering formation to non-Keplerian orbital scenarios. Based on the nonlinear model for relative motion in the CR3BP, two methods for designing teardrop hovering formations along the NRHO have been proposed: one for formations with relatively short revisit distances and another for those with relatively long revisit distances. In the first method, a linear model is used to generate an initial guess for the initial relative velocity, which is then corrected to satisfy the revisit conditions in the nonlinear model. The second method employs a predictor–corrector procedure to continue the hovering trajectories from a short to a long revisit distance. Simulation results verify the effectiveness of the developed methods. The impulse distribution of the design results further providing reference to the parameter selection of the teardrop hovering formation. By considering the dynamical properties near the NRHO, the near-natural 1:1 teardrop hovering formation is achieved. Furthermore, comparison between results obtained from the linear and nonlinear models strengthens the necessity of using the nonlinear model to design the teardrop hovering formation along the NRHO, as the linear model introduces significant errors in the actual relative motion. In future work, we will extend these results to other multi-body periodic orbits and explore teardrop hovering formations with different revisit periods.

\section*{Acknowledgements}
The third author acknowledges the financial support from the National Natural Science Foundation of China (Grant No. 12372044). The fourth author acknowledges the financial support from the National Natural Science Foundation of China (No. U23B6002).

\appendix
\section{Parameter Settings for Fmincon Command}\label{secA1}
The parameter settings for fmincon command are determined through trial and error to ensure both efficiency and accuracy, as presented in Table \ref{tab4}. The ranges of the optimization variables during the optimization are provided in Table \ref{tab5}. Since the minimum revisit distance is set to 1 km, the relative trajectories are considered to satisfy the revisit conditions if the constraint vector $\bm{\psi}$ satisfies $\left|\left|\bm{\psi}\right|\right|<1\times 10^{-9}$.
\begin{table}[!htb]
\caption{Parameter settings for simulations}\label{tab4}%
\centering
\renewcommand{\arraystretch}{1.5}
\begin{tabular}{@{}ll@{}}
\hline
Parameter & Value   \\
\hline
TolX    & $1 \times {10^{ - 11}}$       \\
TolFun    & $1 \times {10^{ - 11}}$       \\
TolCon    & $1 \times {10^{ - 11}}$      \\
MaxIter    & $50000$    \\
MaxFunEvals    & $50000$      \\
\hline
\end{tabular}
\end{table}

\begin{table}[h]
\centering
\renewcommand{\arraystretch}{1.5}
\caption{Optimization variables ranges setting in the fmincon command.}\label{tab5}%
\begin{tabular}{@{}ll@{}}
\hline
NLP variables    & Ranges/(LU/TU) \\
\hline
$\delta u\left(t_0\right)$ & $\delta u\left(t_0\right)\pm 1.5$ \\
$\delta v\left(t_0\right)$ & $\delta v\left(t_0\right)\pm 1.5$ \\
$\delta w\left(t_0\right)$ & $\delta w\left(t_0\right)\pm 1.5$ \\
\hline
\end{tabular}
\end{table}

\bibliographystyle{elsarticle-num}


\printcredits

\end{document}